\documentclass[11pt]{article}
\usepackage[latin1]{inputenc}
\usepackage{amsmath,amsthm,amssymb}

\newtheorem{thm}{Theorem}
\numberwithin{thm}{section}

\def\endproof{\hfill $\Box$\par\vskip3mm}

\def\eq#1{(\ref{#1})}
\def\neweq#1{\begin{equation}\label{#1}}
\def\endeq{\end{equation}}

\def\eq#1{(\ref{#1})}

\def\ep{\varepsilon}

\def\phi{\varphi}
\def\RR{\mathbb R}

\def\di{\displaystyle}
\def\ri{\rightarrow}

\def\incep{\left\{\begin{array}{cl} }
 \def\termin{\end{array}\right. }
\def\2af{2^*_\alpha}

\title{\sc Entire positive solutions of the singular
Emden-Fowler equation with nonlinear gradient term}
\author{\sc Teodora-Liliana Dinu\\
\small Department of Mathematics,
``Fra\c tii Buze\c sti" College, 1100 Craiova, Romania\\ \small Email : {\tt
tldinu@gmail.com}}

\date{}

\begin{document}

\maketitle
\begin{abstract} Let $p$ and $q$ be locally H\"older functions in $\RR^N$,
$p>0$ and $q\geq 0$. We study the Emden-Fowler equation $-\Delta u+
q(x)|\nabla u|^a=p(x)u^{-\gamma}$ in $\RR^N$, where $a$ and $\gamma$ are positive
numbers. Our main result establishes that the above equation has a unique positive solutions
decaying to zero at infinity. Our proof is elementary and it
 combines the maximum  principle for elliptic equations
with a theorem of Crandall, Rabinowitz and Tartar. 
\end{abstract}

\section{Introduction and the main result}
Singular semilinear elliptic problems have been
intensively studied in the last decades.
Such problems arise in the study of non-Newtonian
fluids, boundary layer phenomena for
viscous fluids, chemical heterogenous
catalysts or in the theory of heat conduction in
electrically conducting materials. For instance,
problems
of this type characterize some reaction-diffusion
processes
where the unknown $u\geq 0$ is viewed as the density of a reactant
(see, e.g., \cite{a}). In this framework a major place is played by the
Emden-Fowler singular equation
\neweq{ef}-\Delta u=p(x)u^{-\gamma},\qquad x\in\Omega,\endeq
where $\Omega$ is an open set (bounded or unbounded)
 in $\RR^N$ ($N\geq 3$),  $\gamma>0$, and
 $p:\Omega\ri (0,\infty)$
 is a continuous function. For a comprehensive study of the Emden-Fowler equation
we refer to \cite{d,e,fm,jin,ks,s,wong} and the references therein.
If $\Omega$ is bounded,
Lazer and McKenna proved in \cite{3} that \eq{ef}
has a unique positive solution if $p$ is a
smooth positive function.
The existence of entire positive solutions  for
$\gamma \in (0,1)$ and under certain additional
hypotheses has been established in Edelson \cite{e}
and in Kusano-Swanson \cite{ks}. For instance, Edelson
proved the existence of a solution provided that
$$\int_1^\infty r^{N-1+\lambda (N-2)}\max\limits_{|x|=r}p(x)dr<\infty,$$
for some $\lambda\in (0,1)$. This result is generalized
for any $\gamma >0$ via the sub and super solutions method
in Shaker \cite{s} or by other methods in Dalmasso
\cite{d}.
For further results related to singular elliptic equations we also refer to
\cite{cn1,cn,cp,crt,g}.

The purpose of this paper is to extend some of these results in the more general framework
of singular elliptic equations with nonlinear gradient term.
Problems of this type
arise in stochastic control
theory and have been first studied in Lasry and Lions
\cite{ll}. The corresponding
parabolic equation was considered in Quittner \cite{qu}. Elliptic problems
with nonlinear gradient term have been also studied in various contexts
(see, e.g., \cite{bg,gt,mz,mps}).

We study the problem
 \begin{equation}\label{A}
 \left\{\begin{array}{ll}
 \displaystyle& -\Delta u+q(x)|\nabla u|^a=p(x)u^{-\gamma} \quad\mbox{ in }\ \,
\RR^{N} \\
&\di u>0\quad\mbox{in}\ \,
\RR^N \\
&\displaystyle\lim_{|x|\ri\infty} u(x)= 0,
\end{array}\right.
\end{equation}
where $ N\geq 3$, $a>0$ and $\gamma>0$.
We assume throughout this paper
that  $ p $, $q\in C_{\rm loc}^{0,\alpha}(\RR^{N})$,
 $ p> 0$ and $q\geq 0$ in $\RR^{N}.$
Set
\neweq{phi}\Phi(r)=\max\limits_{|x|=r}p(x).
\end{equation}
 We impose no growth hypothesis on $q$
but we suppose that $p$ satisfies the
following decay condition to zero at infinity:
\neweq{pp}   \int\limits_{0}^{\infty} r\Phi(r) dr<\infty.\endeq
In particular, potentials $p(x)$ which behave like $|x|^{-\alpha}$
as $|x|\ri\infty$, with $\alpha >2$, satisfy this assumption.

 Our main result is the following:

\begin{thm}\label{t1}  Under the above hypotheses, problem \eq{A} has
  a unique  classical solution.
  \end{thm}

   \section{Proof of Theorem \ref{t1}}
We first establish the existence of at least one solution of problem \eq{A}.
For this purpose, for any integer $n\geq 1$, we consider the auxiliary
boundary value problem
\begin{equation}\label{An}
 \left\{\begin{array}{ll}
 \displaystyle& -\Delta u+q(x)|\nabla u|^a=p(x)u^{-\gamma} \quad\mbox{ in }\ \,
B_n \\
&\di u>0\quad\mbox{ in }\ \,
B_n \\
& u= 0, \quad\mbox{ on }\ \,
\partial B_n,
\end{array}\right.
\end{equation}
where $B_n:=\{x\in\RR^N;\ |x|<n\}$.
We observe that the function $\underline u=\ep\varphi_1$ is a subsolution of
\eq{An}, provided that $\ep>0$ is sufficiently small, where $\varphi_1>0$ is
the first eigenfunction of $(-\Delta)$ in $H^1_0(B_n)$. In order to find a supersolution
of \eq{An}, we observe that any solution of
\begin{equation}\label{*}
\left\{
\begin{array}{ll}
\di&  -\Delta u=p(x)u^{-\gamma}\quad \mbox{ in }\ \, B_n \\
\di & u>0\quad \mbox{ in }\ \, B_n \\
\di&  u=0\quad\mbox{ on }\ \, \partial B_n
  \end{array}\right.
  \end{equation}
  is a supersolution of \eq{An}.
 But problem \eq{*} has a solution, by
Theorem 1.1 in Crandall, Rabinowitz and Tartar \cite{crt}.
Denote by $u_n$ this solution.
By standard bootstrap arguments (see \cite{gt1}), $u_n\in C^2(B_n)\cap
C(\overline {B_n})$.
Also, by the maximum principle, it follows that
 $ u_n\leq u_{n+1} $ in $B_n.$ Until now we know that there exists
 $u(x):=\lim_{n\ri\infty}u_n(x)\leq +\infty$, for all $x\in\RR^N$.

Next, we establish the existence of a positive smooth function
$ v $ such that $ u_n\leq v $ in $ \RR^N.$
Let $\Phi$ be defined by \eq{phi} and set
$$ w(r):=\displaystyle K - \int\limits_0^r \zeta^{1-N}\int
\limits_0^{\zeta}\sigma^{N-1}\Phi(\sigma) d\sigma d\zeta ,$$
where
\begin{equation}\label{kmare}
 K:=\displaystyle
\int\limits_0^{\infty}\zeta^{1-N}\int\limits_0^{\zeta}
\sigma^{N-1}\Phi(\sigma)\,d\sigma\,d\zeta \qquad\mbox{for any}\ \, r>0 ,
\end{equation}
provided that the  integral is convergent.
 Then
 $ -\Delta w=\Phi(r) $
and
 $  \lim_{r\to \infty}w(r)=0. $

We prove in what follows that $K<+\infty$. An
 integration by parts yields
 \begin{equation}\label{9}
   \begin{array}{lll}
\di\int\limits_0^r \zeta^{1-N}\int\limits_0^{\zeta}\sigma^{N-1}
\Phi(\sigma)\,d\sigma\,d\zeta =(2-N)^{-1}\int\limits_0^r
\frac{d}{d\zeta}\zeta^{2-N}\int\limits_0^{\zeta}
\sigma^{N-1}\Phi(\sigma)\,d\sigma\,d\zeta  = \\
\di (N-2)^{-1}\left(-r^{2-N}\int\limits_0^r\sigma^{N-1}
 \Phi(\sigma)\,d\sigma+\int\limits_0^r \zeta  \Phi(\zeta)\,d\zeta
\right) \, .
 \end{array}    \end{equation}
Next, by L'Hospital's rule,
 $$\begin{array}{ll}
&\di  \lim_{r\to
\infty}\left(-r^{2-N}\int\limits_0^r
\sigma^{N-1}\Phi(\sigma)\,d\sigma+\int\limits_0^r
\zeta \Phi(\zeta)\,d\zeta \right)=\\
&\di\lim_{r\to \infty}\frac{-\int\limits_0^r
\sigma^{N-1}\Phi(\sigma)\,d\sigma+r^{N-2}\int\limits_0^r \zeta
\Phi(\zeta)\,d\zeta}{r^{N-2}}=
\di \lim_{r\to
\infty}\int\limits_0^r \zeta
\Phi(\zeta)\,d\zeta=\int\limits_0^{\infty} \zeta
\Phi(\zeta)\,d\zeta<\infty,
\end{array}
$$
by our assumption \eq{pp}.
Thus we obtain that $ K=(N-2)^{-1}
\int\limits_0^{\infty} \zeta \Phi(\zeta)\,d\zeta<\infty $.
So, by the definition of $w$,
$ w(r)<(N-2)^{-1} \int\limits_0^{\infty} \zeta
\Phi(\zeta)d\zeta$, for any $r>0 .$

Set
$$v(r):=\left[c(2+\gamma)w(r)\right]^{1/(2+\gamma)},$$
where
$$c:=\left[K(2+\gamma)\right]^{1/(1+\gamma)}.$$
In particular, from $ w(r)\to 0 $ as $ r\to \infty $, we deduce that
 $ v(r)\to 0$ as $r\to \infty .$
Since $w$ is a decreasing function, it follows that $v$ decreases, too.
Hence
$$ \int\limits_{0}^{v(r)}{t}^{1+\gamma}dt\leq
\int\limits_{0}^{v(0)}{t}^{1+\gamma}dt=c w(0)=c K=
\int\limits_{0}^{c}{t}^{1+\gamma}dt. $$
It follows that $ v(r)\leq c $ for all $ r>0 $.

On the other hand,
    $$
\nabla w=\frac{1}{c}\, {v}^{1+\gamma}\nabla v \
\mbox{\ and\ }\
\Delta w=\frac{1}{c}{v}^{1+\gamma}\Delta v +
\frac{1}{c}\left(v^{1+\gamma}\right)'|\nabla v|^2 .
    $$
Hence
\begin{equation}\label{nr}
 \Delta v < c{v}^{-1-\gamma}\Delta w=-c \,
{v}^{-1-\gamma}\Phi(r)\leq -{v}^{-\gamma}\Phi(r)\, .
\end{equation}
By \eq{*} and \eq{nr} we
obtain that $ u_n\leq v $ in $B_n$.
Therefore
$$ u_1\leq u_2\leq \cdots \leq u_n\leq u_{n+1}\leq
 \cdots \leq v  ,$$
with $ v$ vanishing at infinity. Now, standard bootstrap arguments (see \cite{gt1})
imply that $u(x):=\lim_{n\ri\infty}u_n(x)$ is well defined and smooth in $\RR^N$.
Moreover, $u$ is a classical solution of problem \eq{A}.

 \medskip
 We justify in what follows the uniqueness of the solution to problem
 \eq{A}.
Suppose that $u$ and $v$ are arbitrary
 solutions of \eq{A}. In order to establish the uniqueness,
 it is enough to show that
 $ u\leq v $ in $\RR^N$.
Arguing by contradiction, it follows that
  $\max_{x\in\RR^N}
 \left(u(x)- v(x) \right)=:M>0 $. Assume that $ u(x_0)-v(x_0)=M $. Then
 $u(x_0)>v(x_0)>0$,
$ \nabla u(x_0)=\nabla v(x_0)$ and $\Delta (u-v)(x_0)\leq 0$. But
$$
\begin{array}{ll}
\di \Delta (u-v)(x_0)=&\di q(x_0)\left[ |\nabla u(x_0)|^a -
|\nabla v(x_0)|^a\right]+p(x_0)\left(
v^{-\gamma}(x_0)-u^{-\gamma}(x_0)\right)=\\
&\di p(x_0)\left(
v^{-\gamma}(x_0)-u^{-\gamma}(x_0)\right)>0,
\end{array}
$$
which is a contradiction. This implies that $ u\leq v $, and so $u=v$
in $\RR^N$.
\endproof


\begin{thebibliography}{99}{\footnotesize
\bibitem{a} R. Aris, {\it The Mathematical Theory of
Diffusion and Reaction in Permeable
Catalysts}, Clarendon Press, Oxford, 1975.

\bibitem {bg} C. Bandle and E. Giarrusso, Boundary
blow-up for semilinear
elliptic equations with nonlinear gradient terms,
{\it Advances in Differential Equations} {\bf 1} (1996), 133-150.

\bibitem{cn1} A. Callegari and A. Nashman, Some
singular nonlinear equations arising in boundary layer theory, {\it
J.~Math. Anal. Appl.} {\bf 64}
(1978), 96-105.

\bibitem{cn} A. Callegari and A. Nashman, A nonlinear
singular boundary value
problem in the theory of pseudoplastic  fluids, {\it
SIAM J. Appl. Math.} {\bf 38}
(1980), 275-281.

\bibitem{cp} M. M. Coclite and G. Palmieri, On a
singular nonlinear Dirichlet problem,
{\it Comm. Partial Differential Equations} {\bf 14}
(1989), 1315-1327.

 \bibitem{crt} M. G. Crandall, P. H. Rabinowitz and L. Tartar,
 On a Dirichlet problem with a singular nonlinearity,
 {\it Comm. Partial Diff. Equations} {\bf 2} (1977), 193-222.

\bibitem{d} R. Dalmasso, Solutions d'\'equations
elliptiques semi-lin\'eaires singuli\`eres,
{\it Ann. Mat. Pura Appl.} {\bf 153} (1988), 191-201.

\bibitem{e} A. Edelson, Entire solutions of singular
elliptic equations, {\it J.~Math. Anal. Appl.}
 {\bf 139} (1989), 523-532.

\bibitem{fm} W. Fulks and L. S. Maybee, A singular nonlinear equation,
{\it Osaka~J. Math.} {\bf 12} (1960), 1-19.

\bibitem{gt1} D. Gilbarg and N. S.  Trudinger, {\it
Elliptic Partial Differential Equations of Second Order},
Classics in Mathematics, Springer, Berlin, 2001.

\bibitem{g} S. M. Gomes, On a singular nonlinear
elliptic problem,
{\it SIAM J. Math. Anal.} {\bf 17} (1986), 1359-1369.

\bibitem{gt} N. Grenon and C. Trombetti, Existence results for a class of nonlinear
elliptic problems with $p$-growth in the gradient,
{\it Nonlinear Analysis} {\bf 52} (2003), 931-942.

\bibitem{jin} Z. Jin, Solutions for a class of singular semilinear elliptic equations,
{\it Nonlinear Analysis, T.M.A.} {\bf 31} (1998), 475-492.

\bibitem{ks} T. Kusano and C. A. Swanson, Entire positive
solutions of singular semilinear elliptic equations,
{\it Japan J.~Math.} {\bf 11} (1985), 145-155.

\bibitem{ll} J. M. Lasry and P.-L. Lions, Nonlinear
elliptic equations
with singular boundary conditions and stochastic
control with state constraints; the model problem,
{\it Math. Annalen} {\bf 283} (1989), 583-630.

\bibitem{3} A. C. Lazer and P. J. McKenna,
 On a singular nonlinear
 elliptic boundary value problem, {\it Proc. Amer.  Math. Soc.}
 {\bf 111} (1991),  721-730.

\bibitem{mz} H. M\^aagli and M. Zribi, Existence and estimates of solutions
for singular nonlinear elliptic problems, {\it J.~Math. Anal. Appl.}
{\bf 263} (2001), 522-542.

 \bibitem{mps} C. Maderna, C. D. Pagani and S. Salsa, Quasilinear elliptic equations
 with quadratic growth in the gradient, {\it J.~Differential Equations} {\bf 97}
 (1992), 54-70.

 \bibitem{qu} P. Quittner, Blow-up for semilinear
parabolic equations
with a gradient term, {\it Math. Meth. Appl. Sci.}
{\bf 14} (1991), 413-417.

  \bibitem{s} A. W. Shaker, On singular semilinear elliptic
  equations, {\it J. Math. Anal. Appl.}  {\bf 173} (1993), 222-228.

\bibitem{wong} J. S. Wong, On the generalized Emden-Fowler equation,
{\it SIAM Review} {\bf 17} (1975), 339-360.
}
\end{thebibliography}
\end{document}